\documentclass[12pt, leqno]{article}
\usepackage{amssymb,amsmath}
\usepackage[pdftex]{color}
\usepackage{latexsym}
\usepackage[pdftex]{graphicx}
\newtheorem{theorem}{Theorem}
\newtheorem{proposition}{Proposition}

\newtheorem{corollary}{Corollary}

\newtheorem{lemma}{Lemma}

\newcommand{\tx}[1]{\mbox{\;{#1}\;}} 

\newcommand{\R}{\mathbb{R}^n}

\newcommand{\Div}{\hbox{div}\:}
\newcommand{\p}{\partial}
\newcommand{\J}{{\cal{J}}}
\newcommand{\I}{{\cal{I}}}
\newcommand{\E}{{\cal{E}}}
\newcommand{\T}{{\cal{T}}}
\numberwithin{equation}{section}
\begin{document}
\title{Domain perturbations for elliptic problems with Robin boundary conditions of opposite sign.}
\pagestyle{myheadings}
\maketitle
\centerline{\scshape Catherine Bandle}
\medskip
{\footnotesize
\centerline{Mathematische Institut, Universit\"at Basel,}
\centerline{Rheinsprung 21, CH-4051 Basel, Switzerland}
} 
\medskip
\centerline{\scshape Alfred Wagner}
\medskip
{\footnotesize
\centerline{Institut f\"ur Mathematik, RWTH Aachen  }
\centerline{Templergraben 55, D-52062 Aachen, Germany}}
\bigskip

\abstract{We consider the energy of the torsion problem with Robin boundary
conditions in the case where the solution is not a minimizer. Its dependence
on the volume of the domain and the surface area of the boundary is discussed. 
In contrast to the case of positive elasticity constants, the ball does not provide
a minimum. For nearly spherical domains and elasticity constants close to zero  
the energy is largest for the ball. This result is true for general domains in the
plane under an additional condition on the first non-trivial Steklov eigenvalue.
For more general elasticity constants the situation is more involved and it
is strongly related to the particular domain perturbation. The methods used
in this paper are the series representation of the solution in terms of
Steklov eigenfunctions, the first and second shape derivatives and an
isoperimetric inequality of Payne and Weinberger \cite{PaWe61} for the torsional
rigidity.}
\bigskip

\bigskip
\section{Introduction}
Let $\Omega \subset \mathbb{R}^n$ be a  bounded smooth domain  and let $\nu$ denote its outer normal. In this paper we study the Poisson problem
\begin{align}\label{torsion}
\Delta  u+1 =0 \tx{in} \Omega, \quad \p_{\nu }  {u}= \alpha  u \tx{on} \p \Omega.
\end{align}
It is the Euler-Lagrange equation corresponding to the energy functional
$$
E(V,\Omega):=\int_{\Omega} |\nabla V|^2\:dx -\alpha \oint_{\partial \Omega} V^2\:ds -2\int_{\Omega}V\:ds.
$$
If $\alpha<0$ there exists a unique solution $u(x)$ which minimizes  the  energy among all functions in $W^{1,2}(\Omega)$.
In this case  Bucur and Giacomini \cite{BuGi14} have shown that among all domains of given volume the ball has the smallest energy. This property is well-known if $u$ satisfies Dirichlet conditions and follows immediately by symmetrization. The presence of Robin boundary conditions requires completely new arguments. 

In this study we are interesting in the case where $\alpha>0$. The motivation came from the eigenvalue problem $\Delta \varphi +\lambda \varphi =0$ in $\Omega$,  $\p_\nu \varphi= \alpha \varphi$ on $\p\Omega$, considered for the first time by Bareket
\cite{Ba77}. She observed that for nearly circular domains of given area the circle has the largest first eigenvalue. Recently this result was extended to higher dimensions for nearly spherical domains by Ferone, Nitsch and Trombetti \cite{ FeNiTr14} cf. also \cite{BaWa2_14}. The question whether or not the ball is optimal for all domains of the same volume remained opened until recently when Freitas and Krejcirik \cite{FrKr14} showed that annuli have for large $\alpha$ a larger eigenvalue than the ball with the same volume.

If $\alpha >0$ Problem \eqref{torsion} is not always solvable. In fact if $\alpha$ 
coincides with an eigenvalue $0=\mu_1<\mu_2\leq \dots$ of the {\sl Steklov problem}
\begin{align}\label{steklov}
\Delta \phi =0 \tx{in} \Omega, \quad \p_\nu \phi= \mu \phi \tx{on} \p \Omega.
\end{align}
then problem \eqref{torsion} has a solution  if and only if the compatibility condition
\begin{align}\label{compatibility}
\int_{\Omega} \phi_i\:dy =0
\end{align}
holds for all eigenfunctions corresponding to $\mu_i=\alpha$. If \eqref{compatibility} is satisfied then \eqref{torsion} is solvable but the  the solution is not unique.
\smallskip

If $\alpha \neq \mu_i$ then there exists a unique solution.  It is a critical point of $E(V,\Omega)$ in $ W^{1,2}(\Omega)$ in the sense that the Fr\'echet derivative vanishes. However in contrast to the case $\alpha<0$ the critical point is not a local extremum but a saddle point.
\medskip

The goal of this paper is to investigate  $E(u,\Omega)$ among all domains with given volume. In contrast to the case where $\alpha$ is positive the ball has in general not the smallest energy. By means of the shape derivative and a result of Serrin \cite{Se71}for overdetermined boundary value problems it can be shown that the ball is the only critical domain. The analysis of the second shape derivative reveals that for nearly spherical domains and for $\alpha$ small enough the energy is larger or smaller that the one of the ball, depending on the perturbation. The most surprising result in this context is that for $\alpha$ close to zero the ball has the largest energy for all domains of given volume. This phenomenon is related in a wider sense to anti-maximum principles \cite{ClPe79}. At the end we use an upper bound for the torsion of Payne and Weinberger \cite{PaWe61} and obtain an isoperimetric inequality for $E(u,\Omega)$ for all domains in the plane.
\medskip

This paper is organized as follows. First we use the Steklov eigenfunctions to derive a series representation of the energy which will be useful to derive global estimates. This is the content of Section 3. Then we discuss the first shape derivative for general domains and the second shape derivative for nearly spherical domains. At the end we prove the optimality of the disc in two-dimensions. 
\section{Preliminaries}
The Steklov eigenvalues and eigenfunctions will play a crucial role in our considerations. If $\p \Omega$ is Lipschitz continuous then they belong to the Sobolev space $W^{1,2}(\Omega)$ and have a trace in $L^2(\partial \Omega)$. The eigenfunctions can be chosen such that
\begin{align}\label{orthog}
\oint_{\p \Omega} \phi_i \phi_j\:dS = \delta_{ij}, \quad  \int_{\Omega} \nabla \phi_i\cdot \nabla \phi_j \:dx= 0 \tx{if} i\neq j \tx{and} \int_{\Omega} |\nabla \phi_i|^2\:dx= \mu_i.
\end{align}
Moreover every harmonic function $h$ in $\Omega$ with a trace in $L^2(\partial \Omega)$ can be expanded in a series of Steklov eigenfunctions which converges in $W^{1,2}(\Omega)$. It should be mentioned that by a result of Mazya \cite{Ma85} in a Lipschitz domain the norms corresponding to the inner products $<u,v>_\Omega = \int_\Omega \nabla u\cdot \nabla v\:dx +\int_\Omega uv\:dx$ and $<u,v>_{\partial \Omega} = \int_\Omega \nabla u\cdot \nabla v\:dx +\oint_{\partial \Omega} uv\:dS$ are equivalent.

In the next lemma we show how to expand harmonic functions into a series of Steklov eigenfunctions.
\begin{lemma}\label{le:harmonic} (i) Suppose that $\alpha\in \mathbb{R}$ does not coincide with a Steklov eigenvalue $\mu_i$.  Let $h$ be the solution of
\begin{align}\label{eq:harmonic}
\Delta h=0 \tx{in} \Omega, \quad \p_{\nu} h = \alpha h + g(x) \tx{on} \p\Omega.
\end{align}
Then $h= \sum_1^\infty h_i \phi_i$ where
$$
h_i= \frac{\oint_{\p \Omega} \phi_i g \: dS}{\mu_i-\alpha}.
$$
This series converges in $W^{1,2}(\Omega)\cap L^2(\p \Omega)$.
\medskip

(ii) Assume $\alpha= \mu_k$ and denote by ${\mathcal{L}}_k$ the linear space generated by the Steklov eigenfunctions belonging to the eigenvalue $\mu_k$.  A solution exists if and only if the compatibility condition
\begin{align}\label{compatibility:h}
\oint_{\p\Omega} g\phi_k\:dS =0.
\end{align}
is satisfied for all $\phi_k\in {\mathcal{L}}_k$ . In this case  \eqref{eq:harmonic} has infinitely many solutions which are expressed as
$$
h= \sum_{i,i\neq k} h_i \phi_i + {\mathcal{L}}_k,
$$
where $h_i$ is given in (i).
\end{lemma}
{\bf Proof}  Because of the completeness of the Steklov eigenfunctions we can write $h=\sum_1^\infty h_i\phi_i$.
Testing \eqref{eq:harmonic} with $\phi_j$ we get by \eqref{orthog}
$$
0= \int_{\Omega} \phi_j \Delta h \:dy = \oint_{\p \Omega}( \phi_j \p_{\nu} h -h\mu_j\phi_j)\:dS=(\alpha -\mu_j)h_j +\oint_{\p\Omega} \phi_j g\:dS.
$$
This proves the first assertion. The convergence follows from  results by \cite{Ma85,Au04}. The second statement is a consequence of the classical theory on inhomogeneous linear problems.\hfill $\square$
\medskip

For the next considerations we decompose the solution $ u$ of \eqref{torsion} into $ u = h + s$ where $s$ is the solution of the Dirichlet problem
\begin{eqnarray}\label{seq}
\Delta s+1 =0 \:\hbox{in}\: \Omega, \quad s=0 \:\hbox{on}\: \p \Omega.
\end{eqnarray}
Then $h$ is a solution of \eqref{eq:harmonic} with $g=  -\p_{\nu}s$.

A straightforward computation shows that
$$
E(  u,\Omega) = -\int_{\Omega}  u\: dx=  -\int_{\Omega} (h+ s)= - \sum_1^\infty h_i \int_{\Omega}\phi_i\:dx
-\int_{\Omega}s \:dx.
$$
Observe that
$$
-\int_{\Omega} s\:dy= \min_{W^{1,2}_0(\Omega)} \int_{\Omega}\left (|\nabla V|^2 -2V\right)\:dy =:T(\Omega).
$$
Moreover we have
$$
-\sum_1^\infty h_i\int_{\Omega} \phi_i\:dx= \sum_1^\infty\frac{\oint_{\p \Omega}\phi_i\p_{\nu} s\:dS}{\mu_i-\alpha}\int_{\Omega}\phi_i\:dx
$$
and
$$
\int_{\Omega} \phi_i\:dx=  -\int_{\Omega} \phi_i \Delta s\:dx= -\oint_{\p \Omega} \phi_i\p_{\nu}s\:dS.
$$
Hence
\begin{align}\label{id:E}
E(  u,\Omega) =T(\Omega) +\sum_1^\infty \frac{\left (\oint_{\p \Omega} \phi_i\p_{\nu}s\:dS\right)^2}{\alpha -\mu_i}.
\end{align}
\begin{theorem} \label{thm:energy}Assume that $\mu_p<\alpha <\mu_{p+1}$. Set
$$
{\mathcal{E}}^+:= \sum_1^p  \frac{\left (\oint_{\p \Omega} \phi_i\p_{\nu}s\:dS\right)^2}{\alpha-\mu_i}\geq 0 \tx{and} {\mathcal{E}}^-:= \sum_{p+1}^\infty \frac{\left (\oint_{\p \Omega} \phi_i\p_{\nu}s\:dS\right)^2}{\alpha-\mu_i}\leq 0. 
$$
Then the following statements hold true.
\begin{enumerate}
\item
$$
E(  u,\Omega) = T(\Omega) + {\mathcal{E}}^+ +{\mathcal{E}}^-.
$$
\item 
Let ${\bf L_p}$ be the linear space generated by $\{\phi_j\}_{j=1}^p$ and let ${\bf L^\infty_p}$ be the orthogonal complement spanned by
$\{\phi_j\}_{j=p+1}^\infty$. Set 
$$
H(v)=\int_{\Omega }|\nabla v|^2 \:dy-\alpha \oint_{\p \Omega} v^2\:dS + 2\oint_{\p \Omega} v\p_\nu s \:dS.
$$
Then
$$
{\mathcal{E}}^+=\max_{\bf L_p}H(v) \tx{and} {\mathcal{E}}^-=\min_{\bf L_p^\infty}H(v)=\min_vH(v), \tx{where} \oint_{\partial \Omega}v\phi_k\:dS=0 \tx{for} i=1,2..p.
$$
\item If $\alpha = \mu_i$ then \eqref{torsion} has a solution if and only if $\int_{\Omega} \phi_i\:dx= -\oint_{\partial \Omega} \phi_i \p_\nu s\:dS=0$ for all eigenfunctions belonging to $\mu_i$.
\end{enumerate}
\end{theorem}
{\bf Proof} The first assertion follows from \eqref{id:E}. 
Replacing $v \in {\bf L_p}$ by its series $\sum_1^p v_i \phi_i$ we find 
$$
H(v)= \sum_i^p v_i^2(\mu_i-\alpha)+2\sum_1^p v_i\oint_{\p \Omega} \phi_i \p_\nu s \:dS.
$$
By assumption $\mu_i-\alpha <0$ which implies that the positive  maximum is achieved for $v_i=-\oint_{\p \Omega} \phi_i \p_\nu s \:dS/(\mu_i-\alpha)$ which is the Fourier coefficient $h_i$. Inserting this expression into $H(v)$ we obtain $\mathcal{E}^+$.The same argument yields the result for ${\mathcal{E}}^-$. This establishes the second
statement and the last assertion is the compatibility condition stated in \eqref{compatibility}. \hfill $\square$
\medskip

{\sl In the sequel $h_i$ stands for the Fourier coefficient of $h$ in the decomposition $ u= s+h$.}
\medskip

{\sc Remark} The series development \eqref{id:E} holds also for negative $\alpha$.  In this case $\mathcal{E}^+ =0$ and therefore $E( u,\Omega) = T(\Omega) +{\mathcal{E}}^-.$
\medskip

{\sc examples}
1. Le $\Omega= B_R$ be the ball of radius $R$ centered at the origin. The Steklov eigenfunctions for the ball are of the form $r^kX_k(\theta)$ where $\theta\in \p B_1$ and $X_k(\theta)$ are the spherical harmonics of degree $k$. The eigenvalues are 
$\mu=\frac{k}{R}$, $k\in \mathbb{N}$, and their multiplicity is $\frac{(k+n-1)!}{k! (n-1)!}$. By the maximum principle for harmonic functions $\phi_1=$const. is the only radial 
eigenfunction. Here $s=-\frac{R^2}{2n}-\frac{r^2}{2n}$ and thus
$$
h_i (\mu_i-\alpha)= -\oint_{\p B_R} \phi_i \p_\nu s \:dS = 0 \quad \forall i>1.
$$
Consequently  \eqref{torsion} has a solution for all $\alpha\neq 0$. It is of the form 
$$
 u=
\begin{cases}
\frac{R^2}{2n}- \frac{R}{\alpha n} -\frac{r^2}{2n} &\tx{if} \alpha \neq \mu_j ,\\
\frac{R^2}{2n}- \frac{R}{\alpha n} -\frac{r^2}{2n} + w&\tx{if} \alpha = \mu_j, 
\end{cases}
$$
where $w$ is any function in the eigenspace of $\mu_j$.
In both cases we get
\begin{align}\label{e:ball}
E(  u,B_R)= T(B_R) + \frac{|B_R|^2}{\alpha |\p B_R|} =|B_R|\left(-\frac{R^2}{n(n+2)} +\frac{R}{\alpha n}\right).
\end{align}

2. Let $\Omega=\{y:r_0<|y|<R\}$ be an annulus and set $r_0= \kappa R$.
Suppose for simplicity that $\Omega \subset \mathbb{R}^n$, $n>2$. The radial solutions of \eqref{torsion} are of the form
$$
 u=-\frac{r^2}{2n}+c_1 +\frac{c_2}{r^{n-2}}.
$$  
The boundary conditions lead to the linear system
\begin{align}\label{sy:linear}
c_1\alpha +c_2\left(\frac{\alpha}{R} + \frac{n-2}{R^{n-1}}\right)&= \frac{\alpha R^2}{2n} -\frac{R}{n},\\
\nonumber c_1\alpha +c_2\left(\frac{\alpha}{\kappa R} -\frac{n-2}{(\kappa R)^{n-1}}\right)&= \frac{\alpha(\kappa R)^2}{2n} +\frac{\kappa R}{n}
\end{align}
This system has a unique solution if the determinant is different from zero. The determinant vanishes if
$$
\alpha =\alpha_1 =0 \tx{and} \alpha=\alpha_2 = \frac{n-2}{R^{n-2}}\left(\frac{\kappa^{1-n}+1}{\kappa^{-1}-1}\right).
$$
The eigenfunctions of the 
Steklov problem in annulus are similar to those for the ball, namely  $(c_1r^k +c_2r^{-k})X_k(\theta)$, $k=1,\cdots$. In addition to $\phi_1=$const. there is a radial eigenfunction $\phi_r=c_1+c_2r^{2-n}$ with $c_1\mu_r + c_2\left(\frac{\mu_r}{R^{2-n}}+\frac{n-2}{R^{n-1}}\right)=0$ and $c_1\mu_r + c_2\left(\frac{\mu_r}{(\kappa R)^{2-n}}-\frac{n-2}{(\kappa R)^{n-1}}\right)=0$. 
Notice that 
$\alpha_1$ and $\alpha_2$ correspond to the Steklov eigenvalues $\mu_1$ and $\mu_r$ of the radial eigenfunctions. For $\kappa \neq 1$ the inhomogeneous linear system \eqref{sy:linear} is not solvable if
 $\alpha = \mu_r$. Hence  the Fourier
coefficient $h_r$ is not defined. From the symmetry of the annulus it follows that $h_k =0$ for all $k\neq 1,r$.  

The same argument as for the ball shows that in an annulus \eqref{torsion} 
\begin{itemize}
\item has a unique solution if $\alpha \neq \mu_i$,
\item  no solution if $\alpha = \mu_r$,
\item  a family of solutions of the form $-\frac{r^2}{2n}+c_1 +\frac{c_2}{r^{n-2}} +w$ where is in the eigenspace of $\mu_i$  if $\alpha = \mu_i$.
\end{itemize}
Therefore by Theorem \ref{thm:energy} we obtain for the annulus
\begin{align}\label{Eannulus}
E(  u,\Omega) = T(\Omega) +\frac{|\Omega|^2}{\alpha |\p \Omega|} +\frac{h_r^2 }{\alpha-\mu_r}\tx{for all} \alpha \neq 0\tx{and} \alpha\neq \mu_r=\frac{n-2}{R^{n-2}}\left(\frac{\kappa^{1-n}+1}{\kappa^{-1}-1}\right). 
\end{align}
 \section{Global estimates}
3.1. General estimates.
\bigskip

From Theorem \ref{thm:energy} we have  for all $\alpha \neq \mu_j$ that $E(u,\Omega)=T(\Omega) + \mathcal{E}^+ + \mathcal{E}^-$. Many estimates are known for $T(\Omega)$ which is related to the torsion. Less known and more difficult to estimate are the expressions $\mathcal{E}^\pm$. We first start with the observation that
 $$
 a_i:= \oint_{\partial \Omega} \phi_i\partial_\nu s\:dS
 $$
 is the Fourier coefficient of $\partial_\nu s$ with respect to the Steklov eigenfunction $\phi_i$. Let us write
 $$
 \partial_\nu s= \underbrace{\sum_1^p a_i \phi_i}_{\partial_\nu s^+} +\underbrace{\sum_{p+1}^\infty a_i \phi_i}_{\partial_\nu s^-}.
 $$
 Furthermore set for short $\|v\|:=\|v\|_{L^2(\partial \Omega)}$. Then $\|\partial_\nu s^+\|^2= \sum_1^p a_i^2$
 and $\|\partial_\nu s^-\|^2=\sum_{p+1}^\infty a_i^2$.
 Under the assumption $0\leq \mu_p<\alpha <\mu_{p+1}$ it follows immediately that
 \begin{align} \label{e:E+E-}
 \alpha^{-1} \|\partial_\nu s^+\|^2 \leq \mathcal{E}^+\leq (\alpha -\mu_p)^{-1} \|\partial_\nu s^+\|^2,\\
 \nonumber (\alpha-\mu_{p+1})^{-1} \|\partial_\nu s^-\|^2 \leq \mathcal{E}^- \leq (\alpha-\mu_m)^{-1} \sum_{p+1}^m a_i^2.
\end{align}
\medskip

{\sc Application} If $\alpha=-c^2$ is negative we have $\mathcal{E}^+=0$  and therefore \ $\mathcal{E}^-\geq \alpha^{-1}\| \partial_\nu s\|^2_{L^2(\partial \Omega)}$. Hence
$$
E( u,\Omega) \geq T(\Omega)-c^{-2} \|\p_\nu s\|^2.
$$
Equality holds for the balls.  From Schwarz symmetrization it follows immediately that $T(\Omega)\geq T(B_R)$  where $B_R$ is the ball with the same volume as $\Omega$. Also $\int_\Omega |\nabla s_\Omega|^2\:dx \leq \int_{B_R}|\nabla s_{B_R}|^2\:dx$. However it is not clear that
$\|\p_\nu s\|_{L^2(\p\Omega)}\leq \|\p_\nu s\|_{L^2(\p B_R)}$ which would prove that the ball has the smallest energy. Pointwise estimates for $|\nabla s|^2$ are well-known in the literature cf. \cite {GiTr83}, \cite{SpeSch76}.
\bigskip

3.2. Let $0<\alpha<\mu_2(\Omega)$.
\bigskip
  
 In this case 
 $$
 \mathcal{E}^+= \alpha^{-1} \left(\oint_{\partial \Omega} \phi_1\partial_\nu s\:dS\right)^2.
 $$
 Since $\phi_1=\frac{1}{\sqrt{|\partial \Omega|}}.$ we find
 $$
 \mathcal{E}^+= \frac{|\Omega|^2}{\alpha |\partial \Omega|}.
 $$
 This together with Theorem \ref{thm:energy} leads to 
\begin{lemma} \label{le:E} Assume $0<\alpha <\mu_2(\Omega)$. Then
$$
E( u,\Omega) \leq T(\Omega) +\frac{|\Omega|^2}{\alpha |\p\Omega|}.
$$
Equality holds for the ball.
\end{lemma}
If $a_i=0$ for $i=1,\cdots r$, like for instance in the annulus, then the estimate above holds for $0<\alpha<\mu_r(\Omega)$.

 An interesting question is to find an isoperimetric upper bound for
 $$
 \mathcal{J}(\Omega):=T(\Omega) + \frac{|\Omega|^2}{\alpha |\p \Omega|}.
 $$
 If the volume $|\Omega|=|B_R|$ is fixed then -as mentioned before- Schwarz symmetrization implies that $T(\Omega) \geq T(B_R)$,
whereas $\frac{|\Omega|^2}{\alpha |\p\Omega|}\leq \frac{|B_R|^2}{\alpha |\p B_R|}$. The question arises which inequality prevails.
\begin{proposition}\label{le:comp1}
Let $\Omega\neq B_R$ be a fixed domain in $\mathbb{R}^n$ such that $|\Omega|=|B_R|$. Then there exists a positive number $\alpha_0(\Omega)>0$ such that 
$$
\mathcal{J}(\Omega)\begin{cases}
<\mathcal{J}(B_R)&\tx{if} \alpha <\alpha_0,\\
>\mathcal{J}(B_R)&\tx{if} \alpha >\alpha_0.
\end{cases}
$$
\end{proposition}
{\bf Proof}  It is well-known that for any domain different from a a ball $T(\Omega)-T(B_R)= \epsilon_0>0$. Define $\alpha_0= \frac{|B_R|^2}{\epsilon_0}\left( \frac{1}{|\partial B_R|}-\frac{1}{|\partial \Omega|}\right)$. Then
the assertion follows.
\hfill $\square$
\medskip

{\sc remarks} 
\begin{enumerate}

\item  A sharper estimate than in Lemma \ref{le:E} can be derived from Theorem \ref{thm:energy} (2). In fact
$$
E( u,\Omega)\leq T(\Omega) +H(V),
$$
where $V$ is any trial function such that $\oint_{\partial \Omega} V\:dS=0$. Observe that if $V$ is a admissible trial function so is $tV$ for any $t\in \mathbb{R}$. Thus
$$
\mathcal{E}^- \leq \min_{\mathbb{R}} H(tV)= -\frac{\left(\oint_{\partial \Omega} V\partial_\nu s \:dS\right)^2}{ \int_\Omega |\nabla V|^2\:dy -\alpha \oint_{\partial \Omega} V^2\:dS}.
$$
Suppose that the origin is the barycenter with respect to $\partial \Omega$, {\sl i.e.} $\oint_{\partial \Omega} x_i\:dS=0$ for $i=1,\cdots n$. Then $x_i$ is admissible for the variational characterization of $\mathcal{E}^-$.  By our assumption $ \int_\Omega |\nabla x_i|^2 \:dy \geq \mu_2 \oint_{\partial \Omega} x_i^2\:dS >\alpha \oint_{\partial \Omega} x_i^2\:dS $.
Consequently
$$
\mathcal{E}^- \leq - \frac{\sum_1^n\left(\int_ \Omega x_i\:dy\right)^2}{n|\Omega| -\alpha \oint_{\partial \Omega} |x|^2|\:dS}.
$$
\item By the Brock-Weinstock inequality \cite {We54}, \cite{Br01}, $\mu_2(\Omega)\leq \mu_2(B_\rho)$ where $B_\rho$ is the ball of the same boundary measure as $\Omega$, {\sl i.e} $|\partial B_\rho|= |\partial \Omega|$. 
Thus if $|\partial \Omega|$ is  large, $\mu_2(\Omega)$ is small.
\end{enumerate}
3.3. Let $\mu_p(\Omega)<\alpha <\mu_{p+1}(\Omega)$
\bigskip

This case is more involved. From Theorem \ref{thm:energy} it follows that
$E( u,\Omega) =T(\Omega) +\mathcal{E}^++\mathcal{E}^-$.  Rough  estimates are obtained from \eqref{e:E+E-}.

Observe that if the Fourier coefficient $a_p=\oint_{\partial \Omega} \phi_p \partial_\nu s \:dS \neq 0$ then $\mathcal{E}^+ $ is positive and becomes arbitrarily large as $\alpha$ tends to $\mu_p$ from above. Similarly if the Fourier coefficient $\oint_{\partial \Omega} \phi_{p+1} \partial_\nu s \:dS \neq 0$ then $\mathcal{E}^-\neq 0$ then $\mathcal{E}^-$ is negative and becomes arbitraryly small if $\alpha$ tends to $\mu_{p+1}$ from below. 
\medskip

{\sc Examples}

1. In a ball $E( u,\Omega;\alpha)$ has only one pole $\alpha =\mu_1=0$. Hence
$\lim_{\alpha\searrow 0}E( u,\Omega;\alpha)= \infty$ and $\lim_{\alpha\nearrow 0}E( u,\Omega;\alpha)= -\infty$.

2. In an annulus $E( u,\Omega;\alpha)$  has two poles $\alpha =\mu_1=0$ and $\alpha=\mu_r$
see \eqref{Eannulus}.

\section{Domain variations}
\subsection{First domain variation}
\subsubsection{General remarks}
Let $\Omega$ be a family of perturbations of the domain $\Omega$ given by
\begin{equation}\label{not1}
\overline{\Omega}_t=\left\{y:y=x+tv(x)+\frac{t^2}{2}w(x)+o(t^2)\::\: 
x\in\overline{\Omega}\right\},
\end{equation}
where  $v$ and $w$ are smooth vector fields $v,w:\overline{\Omega}\to\R$ belonging to $C^{2,\epsilon}(\overline{\Omega}).$
 
We assume that  on $\p \Omega$, $v$ points in the normal direction, i.e. $v=(\nu\cdot v)\nu$.
The parameter $t$ belongs to $(-t_0,t_0)$ where $t_0$ is chosen so small that $y:\Omega \to \Omega_t$ is a diffeomorphism. 
We consider the family of problems
$$
\Delta_y  u(y,t) +1=0 \tx{in} \Omega_t, \quad  \partial_{\nu_t} u(y,t) =\alpha u(y,t) \tx{on} \partial \Omega_t,
$$
where $\nu_t$ is the outer unit normal at $\Omega_t$. For short we set 
\begin{eqnarray*}
\tilde{u}(t):=u(y(x,t),t)\qquad\hbox{for}\:x\in\Omega,\:\:\vert t\vert< t_0.
\end{eqnarray*}
We now map this problem by means of $y(x,t)$ into 
$\Omega$ and obtain after the change of variable $y \to  x$
\begin{align}\label{eq:torsion2}
\p_j \left(A_{ij}(x,t)\p_j\tilde u(t)\right) +J(t)=0 \tx{in} \Omega, \quad \p_{\nu_A}\tilde u(t)=\alpha m(x,t)\tilde u(t)\tx{on}\p\Omega,
\end{align}
where
\begin{eqnarray*}
\partial_i = \frac{\p}{\p x_i},\quad
dy=J(t)dx,\quad
dS_y=m(t)dS_x,\quad
A_{ij}(t):=\frac{\partial x_i}{\partial y_k}\frac{\partial x_j}{\partial y_k}\:J(t),
\quad\partial_{\nu_A} =\nu_iA_{ij} \partial_j.
\end{eqnarray*}
In \cite{BaWa14} it was shown that for small $|t|$
\begin{eqnarray}\label{jform}
J(t)&:=&\tx{det}(I+tD_v+\frac{t^2}{2}D_w)\\\nonumber&=&1+t\:\Div v+\frac{t^2}{2}\left((\Div v)^2-D_v : D_v+\Div\: w\right)+o(t^2).
\end{eqnarray}
Here we used the notation
\begin{eqnarray*}
D_v : D_v:=\partial_{i}v_j\partial_{j}v_{i},
\end{eqnarray*}
where summation over repeated indices is undestood. Furthermore
$$
m(t)=1 + t (n-1)(v\cdot \nu) H+o(t)
$$
where $H$ is the mean curvature of $\p\Omega$ and 
\begin{eqnarray*}
\Div_{\partial\Omega}v=\Div v-\nu\cdot D_v \nu:=\partial_{i}v_{i}-\nu_{j}\partial_{j}v_{i}\nu_{i}.
\end{eqnarray*}
We also showed that
\begin{eqnarray*}
A_{ij}(0)&=&\delta_{ij};\\
\dot{A}_{ij}(0)&=&\Div v\:\delta_{ij}-\partial_{j}v_{i}-\partial_{i}v_{j};\\
\ddot{A}_{ij}(0)&=&\left((\Div v)^2-D_v:D_v\right)\:\delta_{ij} +2\left(\partial_{k}v_i\:\partial_{j}v_k+\partial_{k}v_j\:\partial_{i}v_k\right)\\
&&
+
2\:\partial_{k}v_i\:\partial_{k}v_j
-
2\:\Div v\left(\partial_{j}v_{i}+\partial_{i}v_{j}\right)
+\Div\:w\: \delta_{ij} -\p_iw_j-\p_j w_i.
\end{eqnarray*}
Similarly we can transform the Steklov problem. In terms of the x-coordinates it reads as
\begin{align}\label{eq:Steklov2}
L_A \phi(t) =0 \tx{in} \Omega, \quad \p_{\nu_A} \phi(t) = \mu(t) m(t) \phi(t) \tx{on} \p \Omega, \quad L_A:= \p_j(A_{ij} \p_i).
\end{align}
The next Lemma is well-known, s. for instance \cite[IV, Sec. 3.5]{Ka} or \cite[VI, Sec. 6]{CoHi65}.
\begin{lemma} Suppose that $\mu_p(\Omega)<\alpha <\mu_{p+1}(\Omega)$. Then there exists $t_0>0$ sich that
$\mu_p(\Omega_t)<\alpha  <\mu_{p+1}(\Omega_t)$ for all $t\in(-t_0,t_0)$..
\end{lemma}
{\bf Proof} By the min-max principle
$$
\mu_p(\Omega_t) =\rm{min}_{L_p}\rm{max}_{V\in L_p} \frac{\int_{\Omega}\nabla V\cdot A(t)\nabla V\:dx}{\oint_{\p\Omega}V^2 m(t)\:dS_x},
$$
where $L_p$ is an $n-$dimensional linear space in $W^{1,2}(\Omega)$. Since $|\nabla V|^2 (1-c_1|t|)<\nabla V\cdot A(t) \nabla V \leq |\nabla V|^2 (1+c_1|t|)$
and $V^2(1-c_2 |t|)<m(t)V^2 \leq V^2(1+c_2 |t|)$. From the min-max principle we obtain that $|\mu_p(\Omega)-\mu_p(\Omega_t)| \leq t C$
where $C$ depends on $v$ and $w$. \hfill $\square$.
\medskip

{\sl In the sequel se shall always assume that $\alpha$ does not coincide with an eigenvalue of $\Omega_t$ for all $t\in (-t_0,t_0)$.}

Suppose that $\Omega \in C^{2,\epsilon}$, $A_{ij}(t)\in C^{1,\epsilon}$, $J(t)\in C^{0,\epsilon}$ and $m(t)\in C^{1,\epsilon} $. We also assume that all the data are at least twice continuously differentiable in $t$. Then by Schauder's regularity theory \cite{GiTr83} it follows that $\tilde u(t) \to \tilde u(0)=:u(x)$ in $C^{2,\epsilon'}$, $\epsilon '<\epsilon$. Our assumptions imply that $|\partial \Omega_t| \to |\partial \Omega|$ which is crucial for the convergence of the eigenvalues.
A  general study of domain perturbations for elliptic problems with Robin boundary conditions is carried out by Dancer
and Daners in \cite{DaDa97}. 
\subsubsection{First variation of the energy}
Consider problem \eqref{torsion} in a class of domains $\Omega$ described in \eqref{not1}. As before the solutions of \eqref{torsion} in $\Omega$ will be denoted by $u(x)$. We shall use the abbreviation ${\mathcal{E}}(t)$ for $E(\tilde u(t),\Omega_t)$. Under the conditions stated above the solution of \eqref{torsion} $\tilde u(t)=\tilde u(y,t)$ is continuous and continuously differentiable in 
$t$. 

\medskip

It was shown in \cite{BaWa14} that the first domain variation  $\frac{d}{dt}{\mathcal{E}}(t)\large |_{t=0}$ is given by
\begin{align*}
\dot{{\mathcal{E}}}(0)=\int_{\p \Omega} (v\cdot \nu)\left[\vert\nabla u\vert^2-2u-2\alpha^2 u^2-\alpha(n-1)u^2H\right]\:dS.
\end{align*}
{\sc Example}
\smallskip

If $\Omega=B_R$ then
\begin{align}\label{edotball}
\dot{\tilde{\mathcal{E}}}(0)=\left(\frac{(n+1)R}{\alpha n^2}-\frac{R^2}{n^2}\right)\int_{\p B_R} (v\cdot \nu)\:dS.
\end{align}

This leads to the following
\begin{corollary} \label{co:critical}Let $\Omega_t$ be a family of nearly spherical domains with prescribed volume $|\Omega_t|=|B_R|$.
Then $\dot{{\mathcal{E}}}(0)=0$. 
\end{corollary}
{\bf Proof} From \eqref{jform} it follows that for volume preserving transformations 
\begin{align}\label{volume1}
\oint_{\p \Omega} (v\cdot \nu)\:dS=0. 
\end{align}
This together with \eqref{edotball} establishes the assertion. \hfill $\square$
\newline
\newline
A further consequence of \eqref{edotball} is the local monotonicity property.
\begin{corollary}
 If $0<\alpha R<n+1$ and $|\Omega_t|>|B_R|$ then $\dot{\cal{E}}(0)>0$, otherwise if $\alpha R>n+1$ then $\dot{\cal{E}}(0)<0$.
 \end{corollary}
 {\bf Proof} By our assumption we have $\int_{\p B_R}(v\cdot \nu)\:dS>0$.  The sign of  $\dot{\cal{E}}(0)$ depends in view of \eqref{edotball} on the sign of $(n+1)\alpha R -(\alpha R)^2$.
 \hfill  $\square$
\subsubsection{First variation of $\mathcal{J}(\Omega_t)$}
In the case $0<\alpha<\mu_2(\Omega)$ (see chapter 3.2) the energy $E(\tilde u,\Omega_t)$ is bounded from above by  ${\cal{J}}(\Omega)=T(\Omega) + \frac{|\Omega_t|^2}{\alpha |\p \Omega_t|}$. 
Let ${\cal{S}}(t)=\vert\partial\Omega_t\vert$. If $|\Omega_t|=|\Omega|$ then the first variation is given by
\begin{eqnarray}
\label{var1}\dot{\J}(0)&=&\dot{T}(0)-\frac{\vert\Omega\vert^2}{\alpha\vert\partial\Omega\vert^2}\dot{{\cal{S}}}(0)
\end{eqnarray}
where
\begin{eqnarray}
\label{var2}\dot{T}(0)&=&-\int\limits_{\partial\Omega}Ê\vert\nabla s\vert^2(v\cdot\nu)\:dSÊ\\
\label{var3}\dot{{\cal{S}}}(0)&=&(n-1)Ê\int\limits_{\partial\Omega}Ê(v\cdot\nu)H\:dSÊ
\end{eqnarray}
Thus for all critical domains the solution $s$ of \eqref{seq} solves the additional boundary condition
\begin{eqnarray}\label{overdet}
\frac{\vert\Omega\vert^2}{\alpha\vert\partial\Omega\vert^2}(n-1)H+\vert\nabla s\vert^2=const.\quad\hbox{on}\:\partial\Omega.
\end{eqnarray}
This is a direct consequence of \eqref{volume1}. By Theorem 3 in \cite{Se71} concerning overdetermined boundary value problems, the ball is the only domain for which on $\p\Omega$, $s$ is constant and $|\nabla s|= c(H)$ for a non-increasing function $c$.
Consequently
\begin{lemma}\label{serrin}
For $\alpha>0$ the ball is the only critical domain for the functional
${\cal{J}}(\Omega)$ among all domains of equal volume.
\end{lemma}
\subsection{Second domain variation for nearly spherical domains}
\subsubsection{Second variation for the energy}
Corollary \ref{co:critical} gives rise to the following question: {\sl  is $E(u, B_R)$ a local extremum among 
the family
 $\Omega_t$, $t\in(-t_0,t_0)$,  of perturbed domains with the same volume as $B_R$?} The answer will be obtained from the second variation.
 
Consider the family of nearly spherical domains $\Omega_t:=\{y=x+tv(x)+\frac{t^2}{2}w(x): x\in \overline{B_R}\}$. Let $\tilde{u}(t):=u(y(x),t)$ be the solution of $\Delta u+1=0$ in $\Omega_t$, $\p_\nu u=\alpha u$ on $\p\Omega_t$ transformed onto $\overline{\Omega}$.
If $\tilde{u}(t)$ is differentiable - this is the case when the data are H\"older continuous as described in the previous section and $\alpha \neq \mu_i(\Omega_t)$ for all $t\in (-t_0,t_0)$ - then
\begin{eqnarray*}
\frac{d}{dt}\tilde{u}(t)\vert_{t=0}=u'(x)+v\cdot\nabla u_0,
\end{eqnarray*}
where $u=u_0$ is the solution of \eqref{torsion} in $B_R$.

It was shown in \cite{BaWa14} that  the shape derivative $u'$ solves the inhomogeneous boundary value problem
\begin{eqnarray}
\label{uprtball1}\Delta u'&=&0\qquad\hbox{in}\:B_R\\
\label{uprtball2}\partial_{\nu}u'&=&\alpha u'+\left(\frac{1-\alpha\:R}{n}\right)v
\cdot\nu\qquad\hbox{on}\:\partial B_R.
\end{eqnarray}
Let us assume that such a solution $u'$ exists. This is certainly the case if $\alpha$ does not coincide with a Steklov eigenvalue $\mu_i(B_R)$.

For the next result we consider perturbations which, in addition to the condition \eqref{volume1}, satisfy the volume preservation of the second order, namely
\begin{align}\label{volume2}
\int_{B_R}((\Div v)^2-D_v : D_v+\Div\:w)\:dx=0.
\end{align}
This formula can be simplified if $v$ points into normal direction only. It takes the form
\begin{eqnarray}\label{vol2}
(n-1)
\int\limits_{\partial\Omega}H(v\cdot\nu)^2\:dS+\int\limits_{\partial\Omega}(w\cdot\nu)\:dS=0.
\end{eqnarray}
Set
\begin{eqnarray*}
Q(u'):= \int\limits_{B_R} |\nabla u'|^2\:dx -\alpha\: \int\limits_{\partial B_R} u'^2\:dS.
\end{eqnarray*}
The following formula was derived in \cite{BaWa14}. Remember that for nearly spherical domains $\mathcal{E}(0)= E(u, B_R)$. Moreover if $\alpha \neq \mu_i(B_R)$,  Lemma 2 implies that for $t$ sufficiently small $\alpha$ never coincides with an eigenvalue $\mu_j(\Omega_t)$.
\begin{lemma}\label{secvar}
Assume $\alpha \neq \mu_i(B_R)$ and let the volume preservation conditions \eqref{volume1} and \eqref{volume2} be satisfied. Put ${\cal{S}}(t):=|\p\Omega|$. Then
\begin{eqnarray}\label{Eradial}
\ddot{{\cal{E}}}(0)=-2 Q(u')
+
\frac{2R}{n^2}(1-\alpha R) \int_{\p B_R}(v\cdot \nu)^2\:dS
-
\frac{R^2}{\alpha n^2}\ddot{\cal{S}}(0) .
\end{eqnarray}
For a ball the second variation of the surface area is of the form
\begin{eqnarray*} 
\ddot{\cal{S}}(0) =\oint_{\partial B_R}\left(|\nabla^{*} (v\cdot\nu)|^2 -\frac{(n-1)}{R^2} (v\cdot \nu)^2\right)\:dS,
\end{eqnarray*}
where $\nabla^{*}$ stands for the tangential gradient on $\p B_R$.
\end{lemma}
\subsubsection{Discussion of the sign of $\ddot{{\cal{E}}}(0)$}
We write for short
\begin{eqnarray}\label{feq}
{\cal{F}}:=-2Q(u')
+
\frac{2R}{n^2}(1-\alpha R) \int_{\p B_R}(v\cdot \nu)^2\:dS.
\end{eqnarray}
In order to estimate ${\cal{F}}$  we consider the Steklov eigenvalue problem \eqref{steklov}
An elementary computation yields $\mu_1=0$,  and $\mu_k =\frac{k-1}{R}$ (for $k\geq 2$ and counted without multiplicity). The second eigenvalue $ \mu_2=1/R$ has multiplicity $n$ and its  eigenfunctions are
$\frac{x_1}{R}, \hdots, \frac{x_n}{R}$. 
\medskip

{\sl From now on we shall count the eigenvalues $\mu_i$ with their multiplicity, i.e. $\mu_2=\mu_3=\mu_{n+1}=1/R$ and $\mu_{n+2}=2/R$ etc.} 

Let $\{\phi_i\}_{i\geq 1}$ be system of Steklov eigenfunctions introduced in Section 2. The function $u'$ solves \eqref{eq:harmonic} with $g=\left(\frac{1-\alpha\:R}{n}\right)v
\cdot\nu$. Hence by Lemma \ref{le:harmonic}
\begin{align*}
u'(x)=\sum_{i=1}^\infty c_i \phi_i \quad\tx{and} \quad(v\cdot \nu)=\sum_{i=1}^\infty b_i\phi_i.
\end{align*}
Note that since the first eigenfunction $\phi_1$ is a constant constant.
the condition 
$$
0=\oint_{\p B_R}(v\cdot \nu)\:dS=\oint_{\p B_R} \phi_1 (v\cdot \nu)\:dS
$$
implies that $b_1=0$. From \eqref{uprtball2} we have also $c_1=0$.
The coefficients $b_i$ for  $i\geq 2$ are determined from the boundary value 
problem \eqref{uprtball1}, \eqref{uprtball2}. In fact
\begin{eqnarray}\label{coeff}
b_i= \frac{n\:c_i\:(\mu_i-\alpha)}{1-\alpha\:R}\qquad\hbox{for} \quad i=2,3,\hdots .
\end{eqnarray}
From the orthonormality conditions \eqref{orthog} of the eigenfunctions it follows that
\begin{align*}
Q(u')&=\sum_{i=2}^\infty c_i^2(\mu_i-\alpha).
\end{align*}
Inserting this into \eqref{feq} we get
\begin{align*}
{\cal{F}}= 2\sum_2^\infty c_i^2\:(\mu_i-\alpha)^2\left[\frac{R}{1-\alpha\:R}-\frac{1}{\mu_i-\alpha}\right].
\end{align*}
Since $\mu_2=\cdots=\mu_{n+1}= \frac{1}{R}$ it follows that
\begin{eqnarray}\label{F}
{\cal{F}}&=& 2\sum_{n+2}^\infty c_i^2\:(\mu_i-\alpha)^2\left[\frac{R}{1-\alpha\:R}-\frac{1}{\mu_i-\alpha}\right]\\
\nonumber&=&2\frac{(1-\alpha R)^2}{n^2}\sum_{n+2}^\infty b_i^2\left[\frac{R}{1-\alpha\:R}-\frac{1}{\mu_i-\alpha}\right].
\end{eqnarray}
\smallskip
Next we shall discuss the sign of $\ddot{\cal{S}}(0)$. Observe that
$$
{\cal{R}}[\chi]=\frac{\oint_{\partial B_R}|\nabla^{*} \chi |^2\:dS}{\oint_{\p B_R} \chi^2\:dS}
$$
is the Rayleigh quotient of the Laplace- Beltrami operator on $\p B_R$. Its eigenvalues $\Lambda_{i_k}$ are $k(n-2+k)/R^2$,
$k\in \mathbb{N}^+$. Observe that the multiplicity of this eigenvalue is the same as for the Steklov eigenvalue corresponding to $k/R$. Remember that for volume preserving perturbations of the first order we  have
 $\oint_{\p B_R} (v\cdot \nu)\:dS=0$ and therefore $(v\cdot \nu)$ is orthogonal to the first eigenfunctionwhich is a constant. Thus
$$
{\cal{R}}[(v\cdot \nu)]\geq \frac{n-1}{R^2}.
$$
Equality holds if and only if $(v\cdot \nu)$ belongs to the eigenspace spanned by $\{\frac{x_i}{R}\}_{i=1}^{n}$. This does not occur if we exclude small translations. Consequently $\ddot{\cal{S}}(0)>0$. This is consistent  with the isoperimetric inequality.
\smallskip

If we replace in $\ddot{\cal{S}}(0)$, $(v\cdot \nu)$ by 
$\sum_1^\infty b_i\phi_i$ we obtain
\begin{align}\label{ddS}
\ddot{\cal{S}}(0)= \sum_2^\infty b_i^2(\Lambda_i- \frac{n-1}{R^2}).
\end{align}
From \eqref{F} and \eqref{ddS} we then get
\begin{align}\label{dee}
\ddot{\mathcal{E}}(0) = \sum_{n+2}^\infty\frac{b_i^2}{\alpha n^2}\underbrace{\left\{2(1-\alpha R)^2(\frac{\alpha R}{1-\alpha R}-\frac{\alpha}{\mu_i-\alpha})-R^2\Lambda_i+n-1)\right\}}_{d_i}.
\end{align}
Since the multiplicity of the Steklov eigenvalues and $\Lambda_i$ depending on $k$ is the same we can replace $\mu_i$ by $k_i/R$ for a suitable integer $k_i$ and $\Lambda_i$ by $k_i(k_i+n-2)/R^2$.
Consequently
\begin{align}\label{d}
d_i=\frac{2\xi(1-\xi)(k_i-1)}{k_i-\xi}-k_i(k_i+n-2)+n-1,
\end{align}
where $ \xi:=\alpha R$ and $k_i=2,3,4 ...$. 
\smallskip

Next we shall discuss the sign of $d_i$. Suppose that $k_p<\xi<k_{p+1}$, $k_p\geq 2$. It is easy to see that 
\begin{align*}
d_i<0   \tx{if} k_ i>k_p \tx{and} i\geq n+2 ,
\end{align*}
If  $\xi =k_p+\epsilon$, $(0<\epsilon<1)$ and $k_p\geq 2$,  then by \eqref{d} 
$$
d_p= \frac{2(k_p+\epsilon)(k_p+\epsilon -1)(k_p-1)}{\epsilon} -k_p(k_p+n-2)+n-1.
$$
For given $k_p\geq 2$ and $n$ we can always find $0<\epsilon$ sufficiently small such that $d_p>0$. 
Observe that for $k_p\geq 2$, 
$d_p-k_p(k_p+n-2)+n-1$ is a monotone decreasing function of $\epsilon$. 
 A lower bound is obtained for $\epsilon=1$, namely
$$
d_p > 2k_p^3-k_p^2-nk_p +n-1.
$$
For $n=2,3,4$ this expression is positive. However in general the sign varies.
\smallskip

\noindent If $k_p=0$, i.e. $0<\xi <1$ or $k_p=1$, i.e. $1<\xi<2$, no positive terms appear in the expression
of $\ddot{\mathcal{E}}(0)$.
\smallskip

\noindent The same situation as for $d_p$ holds for $ d_i<d_p$. Since
$d_i >k_i^2 +(4-n)k_i +n-1$ we have $d_i>0$ if $n\leq 4$.

$k_i<k_p$ the sign of $d_i$ depends on $k_i$ and $n$. If $k_i$ is large compared to $n$ it is positive, otherwise negative. 

These observations are summarized in the following

\begin{lemma}\label{thm1}
(i) Let $0<\alpha R<2$, $\alpha \neq \frac{1}{R}$. Then $\ddot{{\cal{E}}}(0)\leq 0$. Equality holds if and only if $b_i =0$ for all $i\geq n+2$.
\smallskip

(ii) If $k_p<\alpha R<k_{p+1}$, $k_p\geq 2$, then $\ddot{{\cal{E}}}(0)\leq 0$ if $b_i=0$ for $i=n+2,..,k_{p-1}$.
\medskip

(iii) Assume $\alpha R=k_p+\epsilon$, $k_p\geq 2$, $\epsilon\in (0,1)$. Then for every $n$ there exists $\epsilon$ sufficiently small such that $\ddot{{\cal{E}}}(0)\geq 0$ for $b_p\neq 0$ and $b_i=0$ for $i=n+2,..,p-1$ and $b_i=0$ for $i>p$.
\medskip

(iv) If $n\leq 4$ and $k_p<\alpha R<k_{p+1}$ then $\ddot{{\cal{E}}}(0)\geq 0$ if $b_i=0$ for all $i>p$ and arbitrary
$b_i$, $i\leq p$, and $\ddot{{\cal{E}}}(0)\leq 0$ if $b_i=0$ for $i\leq p$ and arbitrary $b_i$, $i>p$.
\end{lemma}
{\sc Example}
\medskip

Let $\Omega\subset \mathbb{R}^2$ be the ellipse whose boundary $\partial \Omega$ is given by 
$$
\left\{\frac{R\cos(\theta)}{1+t},(1+t)R\sin(\theta)\right\},
$$ 
where $(r,\theta)$
 are the polar coordinates in the plane. This ellipse has the same area as the circle $B_R$ and can be interpreted as a perturbation described in \eqref{not1}. We have $y= x+t(-x_1,x_2) + \frac{t^2}{2}(x_1,0) +o(t^2)$. The eigenvalues and eigenfunctions of the Steklov eigenvalue problem \eqref{steklov} in $B_R$ are
$$
\mu = \frac{k}{R} \tx{and} \phi=r^k\begin{cases}
a\cos(k\theta) \\
a\sin(k\theta)
\end{cases}
$$
where $a=\frac{1}{\sqrt{\pi R^{2k+1}}}$ is the normalization constant.
We have
$$
(v\cdot \nu)=-R\cos(2\theta)=b_4\phi_3,
$$
and
\begin{eqnarray*}
 \ddot{S}(0)=\oint_{\partial B_R}\left(|\nabla^{*} (v\cdot \nu)|^2 -\frac{(v\cdot \nu)^2}{R^2}\right)\:dS=3\pi R.
\end{eqnarray*}
A straightforward computation 
yields
\begin{align*}
\ddot{\mathcal{E}}(0)= \left[ -\frac{3}{4\alpha} +\frac{R(1-\alpha R)}{2(2-\alpha R)}\right]\oint_{\partial B_R}(v\cdot\nu)^2\:dS.
\end{align*}
with
$\oint_{\partial B_R}(v\cdot\nu)^2\:dS=\pi R^3$.
From this expression it follows immediately that
\begin{align*}
\ddot{\mathcal{E}}(0)
\begin{cases}
>0 \tx{if} \alpha R>2\\
<0 \tx{if} \alpha R<2.
\end{cases}
\end{align*}
This result is in accordance with Lemma \ref{thm1} (i). In this example $b_i=0$ for all $i\neq 4$. The sign of $\ddot{\mathcal{E}}(0)$ depends therefore on $d_4$. It changes sign at $\alpha R=2$
\bigskip

As we have already mentioned $b_1=0$ for all volume preserving perturbations. The coefficients $b_2,..,b_{n+1}$ belong all to the Steklov eigenvalue $\mu_2=\dots=\mu_{n+1}=1/R$ and give no contribution to $\ddot{\mathcal{E}}(0)$. This is due to the fact that on $\partial B_R$
$$
\sum_2^{n+1} b_i\phi_i =\sum_1^n b_{i+1} c x_i= \vec{b}\cdot \nu,
$$
where $\vec{b}$ is a constant vector. The presence of $b_i$ for $i=2,..,n+1$ means that the perturbed domain
$\Omega$ has been shifted by a vector $t\vec{b}$. Notice that such a shift does not affect the higher coefficients
$b_k$, $k\geq n+2$. Obviously it  leaves the energy invariant. There is therefore no loss in generality to assume that 
\begin{align}\label{ barycenter}
b_2=b_3=\dots=b_{n+1}=0.
\end{align}
This condition also implies that $c_2=\dots =c_{n+1}=0$. Hence Problem \eqref{uprtball1}, \eqref{uprtball2} is solvable for $\alpha R=1$. This observation together with \eqref{dee} implies that for perturbations which are not pure translations or rotations the following result holds true.
\begin{theorem} 
\begin{enumerate}
\item Assume $0<\alpha R<1$. Then
$$
\ddot{\mathcal{E}}(0)\leq - \frac{n-0.5}{\alpha n^2} \oint_{\partial B_R} (v\cdot \nu)^2\:dS<0.
$$
\item Assume $1<\alpha R<2$. Then
$$
\ddot{\mathcal{E}}(0)\leq\frac{1}{\alpha n^2} \left(\frac{2\alpha R(1-\alpha R)}{2-\alpha R}-n-1\right)\oint_{\partial B_R} (v\cdot \nu)^2\:dS<0.
$$
\end{enumerate}
In both cases the energy is  maximal for the ball among all nearly spherical domains of given volume.
\end{theorem}
 In general if $\alpha R>2$ the energy $\mathcal{E}(t)$ has a saddle in $t=0$. 
\begin{theorem}
Assume $n=2,3,4$ and $k_p<\alpha R<k_{p+1}$. Let $\mathcal{L}_p$ be the linear space generated by the eigenfunctions $\phi_i$ belonging to the eigenvalues $\mu_i=1/R,..,k_p/R$ and $\mathcal{L}_p^\perp$ be its complement generated by $\phi_i$ belonging to the remaining eigenvalues $\mu_i= k_{p+1}/R,..., \infty$. Then
$$
\ddot{\mathcal{E}}(0)\begin{cases}
>0 \tx{if} (v\cdot \nu) \in \mathcal{L}_p\\
<0 \tx{if} (v\cdot \nu) \in \mathcal{L}_p^\perp.
\end{cases}
$$
 \end{theorem}
\subsubsection{The second variation of ${\J}(\Omega_t)$}
As for $\ddot{\E}(0)$ we can derive a formula for the second volume preserving domain variation for the functional $\J$. Applying the rules of differentiation we get
\begin{eqnarray}
\label{var4}\ddot{\J}(0)&=& Ê\ddot{\T}(0)+\frac{2\vert\Omega\vert^2}{\alpha{\cal{S}}(\partial\Omega)^3}\dot{{\cal{S}}}^2(0)-\frac{\vert\Omega\vert^2}{\alpha\vert\partial\Omega\vert^2}Ê\ddot{{\cal{S}}}(0)
\end{eqnarray}
In analogy to formulas \eqref{var1} - \eqref{var3} and with the help of \eqref{vol2}
we get
\begin{eqnarray}\label{var5}\ddot{\T}(0)
&=Ê
\int\limits_{\partial\Omega}\vert\nabla s\vert^2
\left((n-1)H(v\cdot\nu)^2-(w\cdot\nu)\right)
\:dSÊ\\
\nonumber &+
2\int\limits_{\Omega}\vert\nabla s'\vert^2\:dx +2\int\limits_{\partial\Omega}(v\cdot\nu)^2\partial_{\nu}s\:dS,
\end{eqnarray}
where  the shape derivative $s'$ satisfies
\begin{eqnarray}\label{shapes}
\Delta s'=0\quad\hbox{in}\:\Omega,\qquad\qquad s'=-v\cdot\nabla s=v\cdot\nu\vert\nabla s\vert\quad\hbox{in}\:\partial\Omega.
\end{eqnarray}
Moreover by formula (2.20) in \cite{BaWa14}
\begin{eqnarray}
\label{var6}\ddot {S}(0)&=&\int\limits_{\partial\Omega}\vert\nabla^{*}(v\cdot\nu)\vert^2\:dS-
\int\limits_{\partial\Omega}\left(\vert A\vert^2-(n-1)^2H^2\right)(v\cdot\nu)^2\:dS\\
\nonumber&&+
(n-1)\int\limits_{\partial\Omega}Ê(w\cdot\nu)HÊ\:dS,
\end{eqnarray}
where
\begin{eqnarray*}
\vert A\vert^2=\sum_{i,j}^{n-1}(\partial_{i}^{*}\nu\cdot x_{\xi_k})(\partial_{k}^{*}\nu\cdot x_{\xi_i}).
\end{eqnarray*}
denotes the socond fundamental form of $\partial\Omega$.
\newline
\newline
From Section 4.1.3 we know that the ball is the only critical point of $\J$. For the ball $B_R$ we have
\begin{eqnarray*}
\dot{S}(0)&=&0,\\
\ddot{S}(0)&=&\oint_{\partial B_R}\left(|\nabla^{*} (v\cdot \nu)|^2 -\frac{n-1}{R^2}(v\cdot \nu)^2\right)\:dS\geq 0,
\end{eqnarray*} 
and
\begin{eqnarray*}
s(x)=\frac{1}{2n}\left(R^2-\vert x\vert^2\right).
\end{eqnarray*}
If $R$ is chosen such that $\vert\Omega_t\vert=\vert B_R\vert$ for all $t\in (-t_0.t_0)$  we get
\begin{eqnarray*}
\ddot{\J}(0)&=&\frac{R^2}{n^2}\int\limits_{\partial B_R}
\left((n-1)H(v\cdot\nu)^2-(w\cdot\nu)\right)
\:dSÊ
+
2\int\limits_{B_R}\vert\nabla s'\vert^2\:dx
-
\frac{2R}{n}\int\limits_{\partial B_R}(v\cdot\nu)^2\:dS\\
&&-
\frac{R^2}{\alpha\: n^2}Ê\:\ddot{{\cal{S}}}(0).
\end{eqnarray*}
The volume constraint \eqref{vol2} then implies
\begin{eqnarray*}
\ddot{\J}(0)=-2\frac{R}{n^2}\int\limits_{\partial B_R}(v\cdot\nu)^2
\:dSÊ
+
2\int\limits_{B_R}\vert\nabla s'\vert^2\:dx
-
\frac{R^2}{\alpha\: n^2}Ê\:\ddot{{\cal{S}}}(0).
\end{eqnarray*}
If we use \eqref{shapes} to eliminate $(v \cdot \nu)$ we can write $\ddot{\J}(0)$ as a functional in $s'$ alone.
\begin{eqnarray*}
\ddot{\J}(0)=\I(s')
:=
2\int\limits_{B_R}\vert\nabla s'\vert^2\:dx
-
\frac{2}{R}\int\limits_{\partial B_R}s'^2\:dS
-
\frac{1}{\alpha}\int\limits_{\partial B_R}\left(\vert\nabla^{*} s'\vert^2 -\frac{n-1}{R^2}s'^2\right)\:dS.
\end{eqnarray*}
\subsubsection{Sign of $\ddot{\J}(0)$}
We like to find the sign of $\I$. For the ball it follows from the volume constraint  that $\oint_{\p B_R} 
s' \:dS=0$. Hence 
$$
\int_{B_R} |\nabla s'|^2\:dx \geq \mu_2(B_R) \oint_{\p B_R} s'^2\:dS.
$$
Since $\mu_2=1/R$ we get the lower estimate
\begin{eqnarray}\label{lowbd}
\I(s')\geq-\frac{1}{\alpha}\ddot{S}(0).
\end{eqnarray} 
Keeping in mind that $s'$ is harmonic we get 
\begin{eqnarray}\label{ests'}
\int_{B_R} |\nabla s'|^2 \:dx =\oint_{\p B_R} s'\partial_{\nu} s'\:dS \leq \frac{1}{2R} \oint_{\p B_R} s'^2\:dS + \frac{R}{2} \oint_{\p B_R}(\p_\nu s')^2\:dS.
\end{eqnarray}
Next we multiply $-\Delta s=1$ with $x\cdot\nabla s$ and integrate over $\Omega$. Since $s=0$ on $\partial\Omega$ this gives    
\begin{eqnarray*}
\int\limits_{B_R}\vert\nabla s'\vert^2\:dx=
\frac{R}{n-2}\int\limits_{\partial B_R}\vert\nabla^{*} s'\vert^2\:dS
-
\frac{R}{n-2}\int\limits_{\partial B_R}(\partial_{\nu} s')^2\:dS
\end{eqnarray*}
If we put this together with the estimate \eqref{ests'} we get
\begin{align*}
\int\limits_{B_R}\vert\nabla s'\vert^2\:dx \leq \frac{1}{nR}\oint_{\p B_R} s'^2\:dS +\frac{R}{2}\oint_{\p B_R} |\nabla ^*s'|^2\:dS.
\end{align*}
This results in the following upper bound.
$$
\mathcal{I}(s') \leq \left( \frac{2R}{n} -\frac{1}{\alpha}\right) \ddot{S}(0).
$$
Thus we have proved
\begin{lemma}\label{maximizer}
For $\alpha<\frac{n}{2R}$ the ball is a local maximizer of $\J(\Omega)$ among nearly circular domains of equal volume.
\end{lemma}
\subsection{Optimality of the ball in two dimensions}
From Proposition \ref{le:comp1} it follows that among all domains of given area $A:=\vert B_R\vert$ the functional $\mathcal{J}(\Omega)$ is smaller than the corresponding expression for the circle provided $\alpha <\alpha_0$ where
$$
\alpha_0=\frac{A^2}{T(\Omega)-T(B_R)}\left(|\p B_R|^{-1}-L^{-1}\right) \tx{where} L=|\p \Omega|.
$$
If we replace $T(\Omega)$ by an upper bound $T^{*}$ then 
$$
\alpha_0\geq \frac{A^2}{T^{*}-T(B_R)}\left(|\p B_R|^{-1}-L^{-1}\right).
$$
Observe that 
$$
T(B_R)=-\frac{A^2}{8\pi}\tx{and} |\p B_R|=\sqrt{4\pi A}.
$$
We are interested in estimates for $T(\Omega)$ which depend only on $L$ and $A$.
in \cite{PaWe61} Payne and Weinberger derived by means of the method of parallel lines such an inequality.
\medskip

Let us introduce the notation
\begin{align*}
2\pi \tilde R:=L, \quad A=:\pi(\tilde R^2- \tilde r^2) \tx{and} \tilde r=y\tilde R.
\end{align*}
Then
$$
y^2=1-\frac{4\pi A}{L^2}, \quad L^2=\frac{4 \pi A}{1-y^2}=\frac{4\pi^2 R^2}{1-y^2} \tx{and} \tilde R^2=\frac{A}{\pi(1-y^2)}=\frac{R^2}{1-y^2}.
$$
Payne-Weinberger's inequality says that
$$
T(\Omega)\leq \frac{\pi}{2}\left(\tilde r^4\log \frac{\tilde r}{\tilde R} -\frac{3}{4}\tilde r^4 +\tilde R^2\tilde r^2-\frac{\tilde R^4}{4}\right).
$$
The expression at the right-hand side is the energy corresponding to the boundary value problem $\Delta U+1=0$ in $B_{\tilde R}\setminus B_{\tilde r}$ with $U=0$ on $\p B_{\tilde R}$ and $\p_\nu U=0$ on $\p B_{\tilde r}$. Consequently equality holds for the disc.

This inequality implies that
\begin{align*}
\epsilon_0= T(\Omega)-T(B_R)\leq  \frac{\pi}{4}\tilde R^4y^2\left[1+y^2\log y^2 -y^2\right]=\frac{\pi R^4}{4(1-y^2)^2}y^2\left[1+y^2\log y^2 -y^2\right].
\end{align*}
Moreover
$$
|\p B_R|^{-1}- L^{-1}= \frac{1}{\sqrt{4\pi A}}(1-\frac{\sqrt{4\pi A}}{L})=  \frac{y^2}{\sqrt{4\pi A}\left(1+\frac{\sqrt{4\pi A}}{L}\right)}=\frac{y^2}{2\pi R(1+\sqrt{1-y^2})}.
$$
Collecting all the terms we obtain the estimate
\begin{align}\label{A0}
\alpha_0\geq \frac{2(1-y^2)^2}{R(1+\sqrt{1-y^2})(1+y^2\log y^2 -y^2)}=:\frac{2}{R}g(y^2).
\end{align}
The function $g(t)$ is monotone increasing for $t\in (0,1)$ with $\lim_{t\to 1}g(t)= 2$ and $g(0)=1/2$. The number $y^2$ measures the defect of $\Omega$ with respect to the circle. The estimate \eqref{A0} together with the monotonicity of $g$ implies
\begin{theorem} \label{J}
(i) Let $\Omega\subset \mathbb{R}^2$ be a domain with fixed area $A$ and let $B_R$ be a disc with the same area. Then
$$
\mathcal{J}(\Omega)\leq \mathcal{J}(B_R) \tx{for all} \alpha \leq \frac{2}{R}g(y^2).
$$
(ii) In particular $\mathcal{J}(\Omega)$ is smaller than the corresponding quantity for the disc if $\alpha\leq 1/R$.
\end{theorem}
Observe that the second statement is consistent with Lemma \ref{maximizer}. As a consequence we have
\begin{corollary}\label{last} Under the same assumptions we have $E(\tilde u,\Omega)$ achieves its maximum for the disc provided $\alpha <\mu_2(\Omega)$.
\end{corollary}
{\bf Proof} From Lemma \ref{le:E} and Theorem \ref{J} (ii) it follows that $E(\tilde u,\Omega) <\min\{\frac{1}{R},\mu_2(\Omega)\}$.
Note that by Weinstock's result $\mu_2(\Omega)\leq \mu_2(B_\rho)\leq\frac{1}{R}$, thus $\min\{\mu_2(\Omega), \frac{1}{R}\}=\mu_2(\Omega)$. \hfill $\square$
\medskip

{\bf Open problem} In order to extend Corollary \ref{last} a generalization of Payne-Weinberger's inequality to higher dimensions would be helpful. This inequality is based on estimates for the length of parallel curves which to our knowledge are not available in higher dimensions.
\bigskip

{\bf Acknowledgement} This paper was initiated during a visit at the Newton Institute in Cambridge. Both authors would like to thank this Institute for the excellent working atmosphere. 

\end{document}